\documentclass[a4paper,12pt]{article}
\usepackage{amsmath}
\usepackage{ulem}
\usepackage{calc}
\usepackage{vmargin}
\usepackage{amstext}
\usepackage{graphicx}
\usepackage{amsfonts}
\usepackage{amssymb}

\begin{document}
\begin{center}
{\Large \bf Iterations and groups of formal transformations}
\end{center}

\begin{center}
	 {\bf O. V. Kaptsov  }\\
	 Institute of Computational Modeling, Siberian Branch, \\ 
	 Russian Academy of Sciences, Krasnoyarsk\\ 
	E-mail: kaptsov@icm.krasn.ru
	 
\end{center}

\noindent
{\bf Abstract}

 In this paper, we consider the problem of formal iteration. We construct an
  area preserving mapping which does not have any square root. 
  This leads to a counterexample to  Moser's existence theorem for an interpolation problem.
  We give examples of formal transformation groups such that  the iteration problem has a solution for every element of the groups.
 


\section{Introduction.}

Iterated functions are objects of study in computer science, fractals, dynamical systems and renormalization group physics \cite{Milnor,Takens,Moser,Kuczma}. 
 Here we will consider continuous iterations of mappings.
Let $\mathbb{K}$ denote either the set of real numbers or the set of complex numbers.
Suppose we are given a local diffeomorphism $u$ of a neighborhood of the origin $0\in \mathbb{K}^n$ onto another and leaves $0$ fixed. 
The problem of continuous iteration consists in finding a one-parameter family  of mappings (a flow) $f(t,x)=f^t(x)$
such that
\begin{equation} \label{group}
f^t\circ f^s= f^{t+s},  \qquad f^1=u, \qquad f^0(x)=x \qquad \forall t,s\in \mathbb{R}.
\end{equation}
The iteration problem  was investigated by Koenigs, Lewis,  Baker, Chen, Sternberg and others.
Bibliographical references can be found in \cite{Kuczma,Chen,GRAMCHEV}.

Every smooth flow $f^t$ is defined by  a system of ordinary differential equations   
$$ y^{\prime}=X(y)  $$
with initial condition $y(0)=x$. Thus the iteration problem is equivalent the following question.
Given a  a local diffeomorphism $u$, does there exist a system of ordinary differential equations such that $y(1)=u$?
If the answer to this question is affirmative then we say that the map $u$ is embedded in the flow $f^t$.

The problem is of great interest in the study of the exponential mapping
of infinite-dimensional Lie algebras of vector fields \cite{Omori,Ovsyannikov,Sternberg}.
Let $\exp(tX)$ denote an one-parameter group generated by a vector field $X$, then the map $\exp: X\mapsto  \exp(X)$ is called the exponential map or time-one map. 
Let $G$ be a group of smooth (or formal) maps,
and  we are given the mapping $u\in G$. 
The question which arises is this: under what conditions is there a vector field $X$ such that $u = \exp(X)$?
If such  a vector field $X$ exists, then it is called the logarithm of $u$. 
We will also say that the formal transformation $u$ possesses a logarithm.

Let us denote by $GS_n(\mathbb{K})$  the group of formal power series transformations \cite{Sternberg}.
Lewis \cite{Lewis} proved that if a transformation $u\in GS_n(\mathbb{K})$ satisfies so-called pseudo-incommensurable
condition, then the iteration problem has a formal  power series solution. This Lewis result has been repeatedly proved by different authors  \cite{Chen,GRAMCHEV,Sternberg}.

In this paper we discuss the iteration problem for some subgroups of the group $GS_n(\mathbb{K})$.
It turns out that there are mappings $u$ to which    
the problem does not even have  formal solution,
namely, we give an example of a polynomial mapping
$u: \mathbb{R}^2\longrightarrow \mathbb{R}^2$ preserving the area such that there does not exist a 2-tuple $g=(g_1,g_2)$ of formal power series $g_1, g_2\in\mathbb{R}[[x,y]]$ with   $g\circ g = u$.
This is a counterexample to Moser's statement \cite{Moser} about the existence of a solution to the iteration problem for area-saving mappings.
We present sufficient conditions for the existence of a solution of the iteration problem.
These conditions allow to indicate some groups of formal transformations such that any element of a group possesses a 
logarithm and the corresponding iteration problem has a formal solution.

\section{Examples and condition for the existence of solutions} 

We begin with the  case of a linear mapping  
$$u(x)=Ux ,  \qquad x\in\mathbb{K}^n, $$ 
where $U$ is an invertible matrix. In this case, a solution of the iteration problem has the form
$$u^t(x)= U^tx=e^{t\ln(U)}x  $$
whenever the matrix $\ln(U)$  is correctly defined. When $\mathbb{K}=\mathbb{C}$ 
the matrix $\ln(U)$ exists but in general it is not unique.
If $\mathbb{K}=\mathbb{R}$ and $U$ is positive definite  then  $\ln(U)$ is a real matrix.
Some details of the linear case can be found in \cite{Lewis}.
Sometimes a nonlinear problem (\ref{group}) can be reduced to a linear one.
This is true if an analytical map $u$ is conjugate to a linear map. 
Some of the most known results in this direction are Poincar\'{e} and Siegel-Sternberg  theorems
\cite{Arnold,Siegel,Sternberg2}.

We now consider the groups of formal transformations.
Let $\mathbb{K}[[x]]$ denote the ring of formal power series in indeterminate $x_1,\dots,x_n$ with coefficients in  $\mathbb{K}$. The ring has a maximal ideal $\mathfrak{M}_1$ and a ideal $\mathfrak{M}_2$ consisting of series without constant and linear terms.  
Denote by $\mathfrak{M}_i^n $ ($i=1,2$) the $n$-ary Cartesian product of $\mathfrak{M}_i$.
Obviously  $\mathfrak{M}^n_1$ is a monoid under substitution of series. 
We denote by $GS_n(\mathbb{K})$ the set of all invertible elements of $\mathfrak{M}^n_1$.
We shall call elements of $GS_n(\mathbb{K})$  formal transformations.
It is clear that $GS_n(\mathbb{K})$ is a group. 
As usual, 
the general linear group of degree $n$ over $\mathbb{K}$ is denoted by $GL_n(\mathbb{K})$.

\noindent
{\bf Example 1.} Let us consider the group $GS_1(\mathbb{C})$ and  a polynomial map
  $$u= e^{i\pi/3}z+z^7.$$ 
It is easy to see that there is no a formal power series
$$g=c_1z+c_2z^2+c_3z^3+c_4z^4+...    $$
 such that  
\begin{equation} \label{g*g}
g\circ g=u.
\end{equation}
Actually, comparing coefficients of $z$ in (\ref{g*g}), we have
$$c_1^2 = e^{i\pi/3}.$$
Then comparing coefficients of $z^2,\dots,z^6$ yields $c_2=\cdots=c_6 =0$.
Finally, comparing coefficients of $z^7$, we obtain
$$ c_1c_7(c_1^6+1) = 1.    $$
This is a contradiction, because $c_1^6+1=0$. 
This example shows that there is no one-parameter group passing through the polynomial $e^{i\pi/3}z+z^7$. If such a group $f^t$  exists, then 
$f^{1/2}\circ f^{1/2}=u$.
But it is not possible  as we just proved.
This example shows that  polynomial map $u= e^{i\pi/3}z+z^7\in GS_1(\mathbb{C})$  does not possess a logarithm.

 We remark that such examples have been known for a long time (see, for example \cite{Kuczma,Sternberg}).

\noindent
{\bf Example 2.}
Let $SS_n(\mathbb{K})$ denote the set
$ \{f\in GS_n(\mathbb{K}): det(Df)=1 \}  ,$
where $Df $ is the Jacobian matrix of $f$, 
i.e. $SS_n(\mathbb{K})$ is a group of volume preserving formal transformations.
Consider an area preserving polynomial mapping $v\in SS_2(\mathbb{R})$  given by
$$ \tilde{x}_1= x_1+x_2^{m+1}, \qquad \tilde{x}_2=x_2$$
and the rotation matrix
$$ M=
\left( \begin{array}{cc}
\cos \alpha  & - \sin\alpha \\
\sin\alpha	& \cos\alpha \\
\end{array} \right),$$
where $\alpha= 2\pi/m$ and $m\geq2$ is an even number.
Thus $u = M v$ is  an area preserving mapping.

It is convenient to use the complex variables $z=x_1+ix_2$ and $\bar{z}=x_1-ix_2$.
Then the mapping $u$ has the form
\begin{equation} \label{g*g2}
u(z,\bar{z})= e^{\frac{2i\pi}{m}}\left(z+\left(\frac{z-\bar{z}}{2i}\right)^{m+1}\right) .
\end{equation}
Let us show that there does not exist a formal series 
$$ g(z,\bar{z})= c_{10}z + c_{01} \bar{z}+c_{20}z^2 +c_{11}z\bar{z} +c_{02}\bar{z}^2 + ... $$
satisfying the condition (\ref{g*g}).
We assume that such series exists and try to find his coefficients.  

Collect all terms belonging to $z,\bar{z}$ in (\ref{g*g}). Then we have two equations
\begin{equation} \label{c10}
e^{i\frac{2\pi}{m}}=c^2_{10}+|c_{01}|^2,
\end{equation}
$$c_{01}(c_{10}+\overline{c}_{10})=0. $$
It follows that
$$c_{01}=0 , \qquad c_{10}=\pm \exp(i\pi/m) .$$
Then comparing coefficients of  $z^k\overline{z}^l$ (
$1<k+l<m+1$)  yields equation 
$$c_{kl}(c_{10} + c^k_{10}\overline{c}^l_{10}) = 0.$$
Obviously, the following inequality holds
$$c_{10}+c^k_{10}\overline{c}^l_{10}\neq0$$
whenever $ 1<k+l<m+1$. Thus we have  $c_{kl}=0$.

Finally, we collect all terms belonging to $z^{m+1}$ and obtain equalities 
$$\frac{\exp(2i\pi/m)}{(2i)^{m+1}}=c_{(m+1)0}c_{10}(1+c^m_{10})=0 ,$$
since $c_{10}=\pm \exp(i\pi/m)$ and $m$ is an even number. 
 This  contradiction proves our assertion.

This example implies that Moser's theorem \cite{Moser} on
the solvability of the iteration problem in the class of formal series is not true even for polynomial mappings. Moreover, it is impossible to find the square root of a  area preserving mapping in the general case. 
This example shows that the polynomial map (\ref{g*g2})  does not possess a logarithm.
We shall see that the above examples are related to resonances.

  

Let $ \lambda_1, \dots, \lambda_n\in\mathbb{C}$ be characteristic values of a matrix 
	$U\in GL_n(\mathbb{K}) $. We recall that an identity of the form 
	\begin{equation} \label{reson}
	\lambda_s=\lambda^{m_1}_1 \cdots \lambda^{m_n}_n \qquad
	m_i\in\mathbb{N}, \qquad \sum_{i=1}^n m_i>1 
	\end{equation}
	is called the resonance (induced by $U$).
We say that the resonance (\ref{reson})
 is not obstructive  if 
\begin{equation} \label{reson^t}
	\lambda_s^t=\lambda^{tm_1}_1 \cdots \lambda^{tm_n}_n  \qquad \forall t\in\mathbb{R}.
		\end{equation}
It is easy to see that we have resonances of the form
$$   \lambda = \lambda^{m+1}  , $$
in Examples 1 and 2 above.
These resonances are obstructive since
$$   \lambda^{\frac{1}{2}} \neq \lambda^{\frac{m+1}{2}}   . $$

Using the  theory of normal forms  we proved  the following statement in \cite{Kap2}.\\
\noindent
{\bf Lemma.} {\it Let $u=Ux+g\in GS_n(\mathbb{C})$ be a formal transformation  with $U\in GL_n(\mathbb{C})$ and $g\in \mathfrak{M}^n_2 $. If any resonance induced by the matrix $U$ is not obstructive then $u$ possesses a logarithm.}

Now we show that the conditions (\ref{reson}), (\ref{reson^t}) are equivalent to Lewis's ones.
Indeed, it follows from (\ref{reson})  that 
$$ \exp(\log\lambda_s) =  \exp(m_1\log\lambda_1+...+m_n\log\lambda_n) . $$
The last equality is equivalent to 
\begin{equation} \label{Lewis1}
   \log(\lambda_s)- \sum_{j=1}^{n} m_j \log(\lambda_j) \in 2\pi i\mathbb{Z}. 	\end{equation}
Similarly, the condition  (\ref{reson^t}) yields 
$$    t( \log(\lambda_s)- \sum_{j=1}^{n} m_j \log(\lambda_j)) \in 2\pi i\mathbb{Z} \qquad \forall t\in \mathbb{R} .$$
It follows that 
\begin{equation} \label{Lewis2}
\log(\lambda_s)= \sum_{j=1}^{n} m_j \log(\lambda_j) . 	
\end{equation}
Conversely, it is easy to see that the equality (\ref{Lewis2}) gives  (\ref{reson^t}) and (\ref{Lewis1}) implies (\ref{reson}).

We recall that Lewis's condition means that any relation (\ref{Lewis1}) implies the equality (\ref{Lewis2}) (see \cite{Sternberg,Lewis}).

One can apply Lemma to obtain subgroups $G$ of $GS_n(\mathbb{K})$ such that any $u\in G$ possesses a logarithm.
For example, consider subgroup $B_l$  which consists of formal transformations 
$$  u=Ux+g ,   \qquad g\in \mathfrak{M}^n_2 , $$
where $U$ is a  lower triangular matrix with real positive eigenvalues.

\noindent
{\bf Corollary.} {\it Any formal transformation $u\in B_l$  possesses a logarithm.}

\noindent
The analogous result holds for subgroup of formal transformations $B^u$ with upper triangular matrices.

{\bf Acknowledgment.} This work is supported by the Krasnoyarsk Mathematical Center and financed by the Ministry of Science and Higher Education of the Russian Federation in the framework of the establishment and development of regional Centers for Mathematics Research and Education (Agreement No. 075-02-2020-1631).


\begin{thebibliography}{39}
	
	\bibitem{Milnor} 
	J. Milnor. Dynamics in one complex variable. Third Edition. Princeton University Press. 2006
	

	\bibitem{Takens}
	 Handbook of dynamical systems, Volume 3. Editors: H. Broer F. Takens B. Hasselblatt. 2010
	
	\bibitem{Moser}
	J. Moser. Lectures on Hamiltonian systems.	Memoirs of the American Mathematical Society,
	no. 81, 1968.
	
	
	\bibitem{Kuczma}
	
	M. Kuczma,  B. Choczewski, R. Ger. Iterative Functional Equations. Cambridge University Press, 1990. 
	
	
	\bibitem{Chen} K. T. Chen, Local Diffeomorphisms--$C^{\infty}$ Realization of Formal Properties // American Journal of Mathematics, Vol. 87, No. 1 (Jan., 1965), pp. 140-157
	
	\bibitem{GRAMCHEV} T. Gramchev, and S. Walcher. Normal Forms of Maps: Formal and Algebraic Aspects // Acta Applicandae Mathematicae (2005) Vol. 87, pp. 123-146
	
	
	
	
	
	\bibitem{Ovsyannikov} L.V. Ovsyannikov. Analytical groups. Novosibirsk, Institute of Hydrodynamics. USSR, 1972
	
	\bibitem{Omori} 
	H. Omori, Infinite-Dimensional Lie Groups, vol. 158 of Translations of Mathematical Monographs, American Mathematical Society, Providence, RI, USA, 1997
	
	

	\bibitem{Sternberg} 
	S. Sternberg. Infinite Lie groups and the formal aspects of dynamic systems// Journal
	of Mathematics and Mechanics, vol. 10 (1961), pp. 451-474.
	
	
	
	
	
	
	\bibitem{Lewis} D. C. Lewis, On formal power series transformations// Duke Mathematical Journal,
	vol. 5 (1939), pp. 794-805.
	
		\bibitem{Arnold} V.I. Arnold, Geometrical Methods in the Theory of Ordinary Differential Equations. 2nd ed. New York etc., Springer-Verlag 1988.
	
	\bibitem{Siegel} C.L. Siegel. Iteration of analytic functions. Ann. Math.(1942) 43, 607-612 .
	
	\bibitem{Sternberg2} S. Sternberg. On the Structure of Local Homeomorphisms of Euclidean n-Space, II//
	American Journal of Mathematics, Vol. 80, No. 3, 1958, pp. 623-631	
	
		\bibitem{Kap2} O. V. Kaptsov. A formal analog of iteration problem//Continuum Mechanics (Dynamica Sploshnoi Sredy in Russin) 1983, Vol 63 , pp. 129-135
	
	
	
	
	
	

	

	
	
\end{thebibliography}
\end{document}